\documentclass[12pt]{article}
\usepackage{times}
\usepackage{booktabs}
\usepackage{pifont}
\usepackage{floatrow}
\usepackage{caption}
\usepackage{mathrsfs}
\usepackage{amsmath}
\usepackage{amsfonts,amsthm,amssymb,mathrsfs,bbding}
\usepackage{txfonts}
\usepackage{tikz}
\usetikzlibrary{positioning}
\usepackage{enumitem}
\usepackage{amsmath}
\usepackage{graphics,multicol}
\usepackage{graphicx}
\usepackage{color}
\usepackage{amssymb}
\usepackage{caption}
\usepackage{cite}
\usepackage{latexsym,bm}
\usepackage{indentfirst}
\usepackage{color}
\usepackage[colorlinks=true,anchorcolor=blue,filecolor=blue,linkcolor=blue,urlcolor=blue,citecolor=blue]{hyperref}
\usepackage{extarrows}
\usepackage{cite}
\usepackage{latexsym,bm}
\usepackage{mathtools}
\evensidemargin=\oddsidemargin\topmargin=-1.5cm

\newtheorem{thm}{Theorem}[section]

\newtheorem{lem}{Lemma}[section]

\theoremstyle{definition}

\addtocounter{section}{0}

\begin{document}
\title{Signless Laplacian spectral conditions for even factors in graphs \footnote{This research was partially supported by National Natural Science Foundation of China
{(No. 12371347, 12501491), the Natural Science Foundation of Hubei Province (Grant
No. 2025AFD006) and the Foundation of Hubei Provincial Department of Education (Grant No. Q20232505).}}}
\author{{\bf Lu Li$^{a}$}, {\bf Hechao Liu$^{a}$}\thanks{Corresponding author.
E-mail addresses: 1851590197@qq.com(L. Li), hechaoliu@yeah.net(H. Liu), hongbo\_hua@163.com(H. Hua), duzn@sxu.edu.cn(Z. Du).}, {\bf Hongbo Hua$^{b}$}, {\bf Zenan Du$^{c}$}\\
{\footnotesize $^a$ Huangshi Key Laboratory of Metaverse and Virtual Simulation, School of Mathematics and Statistics,}\\ {\footnotesize Hubei Normal University, Huangshi, Hubei 435002, P.R. China} \\
{\footnotesize $^b$ Faculty of Mathematics and Physics, Huaiyin Institute of Technology, Huaian, Jiangsu 223003, P.R. China}\\
{\footnotesize $^c$ School of Mathematics and Statistics, Shanxi University, Shanxi, 030006, P.R. China}
}
\date{}

\date{}
\maketitle
{\flushleft\large\bf Abstract}
A spanning subgraph $F$ of a graph $G$ is defined as an even factor of $G$, if the degree $d_F(v)=2k, k\in\mathbb{N}^+$ for every vertex $v\in V(G)$. This note establishes a sufficient condition to ensure that a connected graph $G$ of even order with the minimum degree $\delta$ contains an even factor based on the signless Laplacian spectral radius.

\begin{flushleft}
\end{flushleft}

\section{Introduction}

Let $G=(V(G), E(G))$ be a undirected, finite and simple graph defined by the vertex set $V(G)$ and the edge set $E(G)$. The order and size of a graph $G$ are denoted by $|V(G)|=n$ and $|E(G)|=e(G)$, respectively. For a vertex $v\in V(G)$, $d_G(v)$ represents its degree in $G$, and the minimum degree in $G$ is given by $\delta(G)=\min\{d_G(v): v\in V(G)\}$. The number of odd components in $G$ is written as $o(G)$. For a vertex subset $S\subseteq V(G)$, the subgraph induced by $S$ is denoted by $G[S]$, while the subgraph induced by $V(G)\setminus S$ is denoted by $G-S$. Furthermore, $K_n$ denotes the complete graph of order $n$. For two vertex-disjoint graphs $G_1$ and $G_2$, their graph union is denoted by $G_1\cup G_2$, and their join is denoted by $G_1\vee G_2$, is the graph formed by adding every possible edge between the vertex sets $V(G_1)$ and $V(G_2)$ to $G_1\cup G_2$.

A spanning subgraph $F$ of a graph $G$ is defined as an even factor of $G$, if the degree $d_F(v)=2k, k\in\mathbb{N}^+$ for every vertex $v\in V(G)$. A special case of an even factor is a 2-factor, which occurs when every vertex in the spanning subgraph has a degree of exactly 2, that is, $d_F(v)=2$ for all $v\in V(G)$.

A cornerstone in the theory of graph factors is Tutte's 1-factor theorem \cite{W1947The}, established in 1947. And Tutte \cite{W1954A} also established a characterization for the existence of a 2-factor in a graph. The foundational result has inspired extensive research and spurred the development of numerous related findings in the study of graph factors. Ota and Tokuda \cite{K1996A} demonstrated that any star-free ($K_{1,n}$-free) graph $G$ with minimum degree $\delta(G)\geq2n-2$ exists a 2-factor for integer $n\geq 3$. Ryj\'{a}\v{c}ek, Saito and Schelp \cite{Z1999Closure} later showed that a claw-free ($K_{1,3}$-free) graph $G$ contains a 2-factor if and only if its closure does. Steffen and Wolf \cite{E2022Even} proved that any $k$-critical graph contains an even factor if it has at most $2k-6$ vertices of degree 2. Further sufficient conditions were provided by Chen and Chen \cite{X2023Local} for connected graphs of order at least 3 to contain a 2-factor. Xiong \cite{L2017Characterization} established two necessary and sufficient conditions for a specific class of graphs to contain an even factor. Additional results on even factors in graphs have been provided by Kobayashi and Takazawa \cite{Y2009Even}, Zhang and Xiong \cite{L2025Characterizing}.

Let $G$ be a graph with vertex set $V(G)=\{v_1, v_2, ...,v_n\}$. Its adjacency matrix $A(G)$ is a $n\times n$ symmetric matrix with entries $a_{ij}$, where $a_{ij}=1$ if $v_i$ is adjacent to $v_j$, and $a_{ij}=0$ otherwise. The spectral radius of $G$, denoted by $\rho(G)$, is defined as the largest eigenvalue of its adjacency matrix $A(G)$. The signless Laplacian matrix $Q(G)$ is defined as $Q(G)=D(G)+A(G)$, where $D(G)$ is the diagonal degree matrix of $G$. The signless Laplacian spectral radius of $G$, denoted by $q(G)$, is defined as the largest eigenvalue of its signless Laplacian matrix $Q(G)$.

O \cite{S2021Spectral} established a lower bound on the spectral radius (resp. size) of a connected graph $G$ that guarantees the presence of a perfect matching. Zhou, Bian and Wu \cite{SSufficient} established a spectral radius (resp. size) condition to ensure that a connected graph $G$ of even order with the minimum degree $\delta$ contains an even factor. Zhou \cite{zhouSufficient} provided a spectral radius (resp. size) condition for a graph with minimum degree to be $k$-critical with respect to $[1,b]$-odd factor. Zhou and Zhang \cite{S2025Signless} established a lower bound on the signless Laplacian spectral radius of $G$ to ensure that $G$ is $k$-extendable. Wang, Yang and Yang \cite{T2025Signless} provided a signless Laplacian spectral radius condition to guarantee that $G$ is $k$-critical with respect to $[1,b]$-odd factor. For additional findings concerning the relationship between spectral radius and spanning subgraphs, we refer readers to references \cite{J2024Spectral, Y2025Binding, S2025Some, S2025A, L2024Some}.

Motivated by \cite{SSufficient, zhouSufficient, T2025Signless} directly, this work establishes a sufficient condition in terms of signless Laplacian spectral radius to guarantee that a connected graph $G$ of even order with the minimum degree $\delta$ contains an even factor. We have the following result:
\begin{thm}
Let $G$ be a connected graph of even order $n\geq 7\delta-7$ with minimum degree $\delta\geq 2$. If
\begin{center}
$q(G)\geq q(K_\delta \vee (K_{n-2\delta +1} \cup (\delta -1)K_1))$,
\end{center}
then $G$ contains an even factor, unless $G=K_\delta \vee (K_{n-2\delta +1} \cup (\delta -1)K_1)$.
\end{thm}

\section{Preliminaries}
This section presents several preliminary lemmas that are instrumental in proving our main theorem. Among these, a key result by Yan and Kano \cite{Z2020Strong} provides a sufficient condition for the existence of even factors in graphs.

\begin{lem} \label{lem01} (\cite{Z2020Strong})
If $G$ is a graph of even order $n$, then $G$ contains an even factor if and only if
\begin{center}
$o(G-S)< |S|$
\end{center}
holds for any $S\subseteq V(G)$ with $|S|\geq 2$.
\end{lem}

\begin{lem} \label{lem02} (\cite{L2024Some})
Let $s, p, n_1, n_2, \ldots, n_t$ be positive integers such taht $t\geq 2, 1\leq p\leq n_1\leq n_2\leq \ldots\leq n_t$, and $\sum_{i=1}^{t}n_i=n-s$. Then
\begin{center}
$q(K_s\vee (K_{n_1}\cup K_{n_2}\cup \cdots\cup K_{n_t}))\leq q(K_s \vee ((t -1)K_p\cup K_{n-s-p(t-1)}))$,
\end{center}
with equality if and only if $n_1= n_2= \ldots= n_{t-1}= p$.
\end{lem}

\begin{lem} \label{lem03} (\cite{X2013Matrix})
Let $H$ be a subgraph of a connected graph $G$. Then
\begin{center}
$q(G)\geq q(H)$,
\end{center}
with equality if and only if $G=H$.
\end{lem}

Consider a real matrix $M$ of order $n$, and let $\mathcal{N}=\{1, 2, \cdots, n\}$. Given a partition $\pi :\mathcal{N}=\mathcal{N}_1\cup \mathcal{N}_2\cup \cdots\cup \mathcal{N}_r$, the matrix $M$ can be expressed in block form as
\begin{center}
$\left.M=\left( \begin{array} {cccc}M_{11} & M_{12} & \cdots & M_{1r} \\ M_{21} & M_{22} & \cdots & M_{2r} \\ \vdots & \vdots & \ddots & \vdots \\ M_{r1} & M_{r2} & \cdots & M_{rr} \end{array}\right.\right)$,\\
\end{center}
where each block $M_{ij}$ formed by rows in $\mathcal{N}_i$ and columns in $\mathcal{N}_j$. Let $b_{ij}$ denote the average row sum of $M_{ij}$, obtained by dividing the sum of all its entries by the number of rows in $M_{ij}$. The matrix $M_\pi=(b_{ij})_{r\times r}$ is referred to as the quotient matrix of $M$ with respect to the partition $\pi$. If every row within each block $M_{ij}$ has an identical sum, the partition $\pi$ is termed equitable. When this condition is met, the corresponding quotient matrix $M_\pi$ is specifically referred to as the equitable quotient matrix.

The following lemmas establish fundamental spectral properties of matrices with equitable partitions and their applications to graph structures.

\begin{lem} \label{lem04} (\cite{L2019On})
If $M$ is a real matrix of order $n$ with an equitable partition $\pi$, and $M_\pi$ is the corresponding equitable quotient matrix, then every eigenvalue of $M_\pi$ is also an eigenvalue of $M$. Furthermore, if $M$ is nonnegative, the largest eigenvalues of $M$ and $M_\pi$ coincide.
\end{lem}

\begin{lem} \label{lem05} (\cite{T2025Signless})
If $M$ is a nonnegative irreducible matrix of order $n$ with an equitable partition $\pi$, and $X$ is the Perron vector of $M$, then the entries of $X$ are constant on each cell of the partition $\pi$.
\end{lem}

\begin{lem} \label{lem06} (\cite{JSpectral})
Let $G=K_s\vee (K_{n_1}\cup K_{n_2}\cup \cdots\cup K_{n_t})$, where $t\geq 2$ and $n_1\leq n_2\leq \ldots\leq n_t$. Then $G$ exists a natural equitable vertex partition $\pi=\{V(K_s), V(K_{n_1}), V(K_{n_2}), \ldots, V(K_{n_t})\}$. Let $X$ be the Perron vector of the matrix $\alpha D(G)+A(G)$ for $\alpha=0,1$, with $x_i$ denoting the common value for vertices in $K_{n_i} (1 \leq i \leq t)$. Then we have $x_i \leq x_{i+1}$, with equality if and only if $n_i = n_{i+1}$.
\end{lem}

\section{The proof of Theorem 1.1}

In this section, we present the proof of Theorem 1.1, which elucidates the connection between the signless Laplacian spectral radius and the existence of even factors in graphs.
\begin{proof}
Suppose, by contradiction, that $G$ does not contain an even factor. Then by Lemma \ref{lem01}, there exists a nonempty subset $S\subseteq V(G)$ with $|S|\geq 2$ such that $o(G-S)\geq |S|$. Let $s=|S|$. It follows that $G$ is a spanning subgraph of the graph $G_1 = K_s\vee (K_{n_1}\cup K_{n_2}\cup \ldots\cup K_{n_s})$, where each $n_i$ is a positive odd integer satisfying $1\leq n_1\leq n_2\leq \ldots\leq n_s$ and $\sum_{i=1}^{s}n_i=n-s$. By Lemma \ref{lem03}, we obtain the inequality
\begin{equation}\tag{3.1}
q(G)\leq q(G_1),
\end{equation}
with equality if and only if $G = G_1$. The remainder of the proof proceeds by analyzing three distinct cases according to the value of $s$.

\textbf{Case 1.} $s\geq \delta+1.$

We define $G_2=K_s \vee (K_{n-2s +1} \cup (s -1)K_1)$, where $n\geq 2s$. By Lemma \ref{lem02}, we obtain the inequality
\begin{equation}\tag{3.2}
q(G_1)\leq q(G_2),
\end{equation}
with equality if and only if $G_1=G_2$. Under the vertex partition $V(G_2)=V(K_s)\cup V(K_{n-2s+1})\cup V((s-1)K_1)$, the quotient matrix of $Q(G_2)$ is given by
\begin{center}
$B_2= \begin{pmatrix} n+s-2 & n-2s+1 & s-1 \\ s & 2n-3s & 0 \\ s & 0 & s \end{pmatrix}$ ,
\end{center}
and its characteristic polynomial is
\begin{center}
$\varphi_{B_2}(x)=x^3-(3n-s-2)x^2+(2n^2+ns-4n-4s^2+4s)x-2sn^2+4ns^2+2ns-2s^3-2s^2$.
\end{center}
Note that the partition $V(G_2)=V(K_s)\cup V(K_{n-2s+1})\cup V((s-1)K_1)$ is equitable. By Lemma \ref{lem04}, the largest root of $\varphi_{B_2}(x)=0$ equals $\rho(G_2)$.

We define $G_*=K_\delta \vee (K_{n-2\delta +1} \cup (\delta -1)K_1)$. Under the vertex partition $V(G_*)=V(K_\delta)\cup V(K_{n-2\delta+1})\cup V((\delta-1)K_1)$, the quotient matrix of $Q(G_*)$ is given by
\begin{center}
$B_*= \begin{pmatrix} n+\delta-2 & n-2\delta+1 & \delta-1 \\ \delta & 2n-3\delta & 0 \\ \delta & 0 & \delta \end{pmatrix}$ ,
\end{center}
and its characteristic polynomial is
\begin{center}
$\varphi_{B_*}(x)=x^3-(3n-\delta-2)x^2+(2n^2+n\delta-4n-4\delta^2+4\delta)x-2\delta n^2+4n\delta^2+2n\delta-2\delta^3-2\delta^2$.
\end{center}
Note that the partition $V(G_*)=V(K_\delta)\cup V(K_{n-2\delta+1})\cup V((\delta-1)K_1)$ is equitable. By \ref{lem04}, the largest root of $\varphi_{B_*}(x)=0$ equals $\rho(G_*)$.
Through direct computation, we obtain the polynomial difference
\begin{equation}\tag{3.3}
\varphi_{B_2}(x)-\varphi_{B_*}(x)=(s-\delta)f(x),
\end{equation}
where $f(x)=x^2+(n-4s-4\delta+4)x-2n^2+2n(2s+2\delta+1)-2(s^2+s\delta+\delta^2)-2(s+\delta)$.

Since $K_{n-\delta+1}$ is a proper subgraph of $G_*$, it follows from \ref{lem02} that
\begin{center}
$q(G_*)> q(K_{n-\delta+1})= 2n-2\delta$.
\end{center}
Next we claim $\varphi_{B_2}(x)>\varphi_{B_*}(x)$ for all $x\in[2n-2\delta, +\infty]$.

The axis of symmetry of the quadratic function $f(x)$ is located at $x=\frac{4s+4\delta-n-4}{2}$. Since $\delta+1\leq s\leq \frac{n}{2}$ and $n\geq 7\delta-7$, we derive the inequality
\begin{center}
$(2n-2\delta)\times2-(4s+4\delta-n-4)= 5n-8\delta-4s+4\geq 3n-8\delta+4\geq 13\delta-17>0$.
\end{center}
This implies $\frac{4s+4\delta-n-4}{2}< 2n-2\delta< q(G_*)$, and it follows that $f(x)$ is monotonically increasing on the interval $[2n-2\delta, +\infty)$. Hence,
\begin{align*}
f(x)&\geq f(2n-2\delta)\\&= 4n^2-(4n-6\delta)s-14n\delta+10n+10\delta^2-10\delta-2s^2-2s\\
&\geq 4n^2-(4n-6\delta)\frac{n}{2}-14n\delta+10n+10\delta^2-10\delta-2(\frac{n}{2})^2-2\frac{n}{2}\\
&= \frac{3}{2}n^2-(11\delta-9)n+10\delta^2-10\delta\\
&\geq \frac{3}{2}(7\delta-7)^2-(11\delta-9)(7\delta-7)+10\delta^2-10\delta\\
&= \frac{1}{2}(13\delta^2-34\delta+21)\\
&\geq \frac{1}{2}(13\cdot2^2-34\cdot2+21)\\
&= \frac{5}{2}>0 \tag{3.4}
\end{align*}
Based on (3.3) and (3.4), together with the condition $s\geq\delta+1$, it follows that $\varphi_{B_2}(x)>\varphi_{B_*}(x)$ for all $x\geq2n-2\delta$. Since $q(G_*)> 2n-2\delta$, we have $\varphi_{B_2}(x)>\varphi_{B_*}(x)\geq0$ for all $x\geq q(G_*)$. Therefore, the equation $\varphi_{B_2}(x)=0$ has no roots in the interval $[q(G_*), +\infty)$, implying $q(G_2)< q(G_*)$. Combining this result with (3.1) and (3.2), we obtain the chain of inequalities
\begin{center}
$q(G)\leq q(G_1)\leq q(G_2)< q(G_*)=q(K_\delta \vee (K_{n-2\delta +1} \cup (\delta -1)K_1))$,
\end{center}
which contradicts the initial assumption that $q(G)\geq q(K_\delta \vee (K_{n-2\delta +1} \cup (\delta -1)K_1))$.

\textbf{Case 2.} $s=\delta$.

Recall that $G_1 = K_s\vee (K_{n_1}\cup K_{n_2}\cup \cdots\cup K_{n_s})$. By Lemma \ref{lem02} and the condition  $s=\delta$, we derive
\begin{center}
$q(G_1) \leq q(K_\delta \vee (K_{n-2\delta +1} \cup (\delta -1)K_1))$,
\end{center}
with equality if and only if $G_1 = K_\delta \vee (K_{n-2\delta +1} \cup (\delta -1)K_1)$. Combining this with inequality (3.1), we conclude
\begin{center}
$q(G) \leq q(K_\delta \vee (K_{n-2\delta +1} \cup (\delta -1)K_1))$,
\end{center}
where equality holds if and only if $G = K_\delta \vee (K_{n-2\delta +1} \cup (\delta -1)K_1)$, again leading to a contradiction.

\textbf{Case 3.} $s\leq \delta-1$.

We define $G_3=K_s \vee (K_{n-s-(\delta+1-s)(s-1)} \cup (s -1)K_{\delta+1-s})$. Recall that $G$ is a spanning subgraph of $G_1 = K_s\vee (K_{n_1}\cup K_{n_2}\cup \ldots\cup K_{n_s})$, where $n_1\leq n_2\leq \ldots\leq n_s$ and $\sum_{i=1}^{s}n_i=n-s$. Obviously, $n_1\geq \delta+1-s$ because $\delta(G_1)\geq \delta(G)=\delta$. In terms of Lemma \ref{lem02}, we conclude
\begin{equation}\tag{3.5}
q(G_1)\leq q(G_3),
\end{equation}
with equality holding if and only if $(n_1, n_2, \ldots, n_s)=(\delta+1-s, \delta+1-s, \ldots, n-s-(\delta+1-s)(s-1))$.
Partition the vertex set of $G_3$ as follows:\\
$\begin{aligned} & S=V(K_{s})=\{u_{1}, u_{2}, \ldots, u_{s}\}, \\ & V_{1}=\bigcup_{i=1}^{s-1}V\left(K_{\delta+1-s}^{(i)}\right)=\{v_{i,j}: 1\leq i\leq s-1, 1\leq j\leq\delta+1-s\}, \\ &
V_{2}=V\left(K_{n-s-(\delta+1-s)(s-1)}\right)=\{w_{1}, w_{2}, \ldots,w_{m}\}, \mathrm{where~}m=n-s-(\delta+1-s)(s-1). \end{aligned}$
\vspace{0.5em}

Let $X$ be the Perron vector of $Q(G_3)$. Lemma \ref{lem05} ensures that the entries of $X$ are constant on each of the vertex subsets $S, V_1$, and $V_2$. Denote these constant values by $x_1, x_2,$ and $x_3$, respectively. Furthermore, by Lemma \ref{lem06}, we have $x_2\leq x_3$.

We construct $G_4$ from $G_3$ by selectively removing edges within the cliques $K_{\delta+1-s}^{(i)}$ and simultaneously adding edges to merge vertices into larger cliques. Formally:\\
\vspace{0.5em}
$\begin{aligned} E_{\mathrm{removed}} &=E\left(K_{\delta+1-s}^{(1)}\right)\cup\bigcup_{i=2}^{s-1}\{v_{i,1}v_{i,j}: 2\leq j\leq\delta+1-s\}, \\ E_{\mathrm{added}}&=E_1\cup E_2, \end{aligned}$\\
\vspace{0.5em}
where\\
\vspace{0.5em}
$\begin{aligned} & E_{1}=\{vw_{r}: v\in V_{1},1\leq r\leq\delta-s\}, \\ & E_{2}=\{v_{i,j}w_{r}: 2\leq i\leq s-1, 2\leq j\leq\delta+1-s,\delta-s+1\leq r\leq m\}. \end{aligned}$\\

Thus, $G_4=G_3+E_{\mathrm{added}}-E_{\mathrm{removed}}$. Next we claim that $q(G_4)>q(G_3)$. The quadratic form difference is:
\begin{center}
$X^T(Q(G_4)-Q(G_3))X=\underbrace{|E_{\mathrm{added}}|(x_2+x_3)^2}_{\mathrm{Gain}}
-\underbrace{|E_{\mathrm{removed}}|(x_2+x_2)^2}_{\mathrm{Loss}}.$
\end{center}
where\\

$\begin{aligned}  |E_{\mathrm{added}}|&=(s-1)(\delta+1-s)(\delta-s)+(m-\delta+s)(s-2)(\delta-s), \\ |E_{\mathrm{removed}}|&=\binom{\delta+1-s}{2}+(s-2)(\delta-s). \end{aligned}$\\
Substituting into the quadratic form:\\

$\begin{aligned} & X^T(Q(G_4)-Q(G_3))X \\ & = [(s-1)(\delta+1-s)(\delta-s)+(m-\delta+s)(s-2)(\delta-s)] (x_2+x_3)^2 \\ & -\left[\binom{\delta+1-s}{2}+(s-2)(\delta-s)\right](x_2+x_2)^2 \\ & =\left[2(s-1)\binom{\delta+1-s}{2}(x_{2}+x_{3})^{2}-\binom{\delta+1-s}{2}(x_{2}+x_{2})^{2}\right] \end{aligned}$

$\begin{aligned} &+\left[(m-\delta+s)(s-2)(\delta-s)(x_2+x_3)^2-(s-2)(\delta-s)(x_2+x_2)^2\right] \\ & \geq (2s-3)\binom{\delta+1-s}{2}(x_{2}+x_{2})^{2}+(m-\delta+s-1)(s-2)(\delta-s)(x_2+x_2)^2. \end{aligned}$

By assumptions, we have $s\geq2, \delta-s\geq1, \binom{\delta+1-s}{2}\geq1$ and $m-\delta+s-1\geq0$. Thus $X^T(Q(G_4)-Q(G_3))X>0$, implying $q(G_4)\geq X^TQ(G_4)X>X^TQ(G_3)X=q(G_3)$. This proves that $q(G_4)>q(G_3)$.

Since $G_4\subseteq G_*$, Lemma \ref{lem03} implies $q(G_4)\leq q(G_*)$. Combining this result with (3.1) and (3.5), we obtain the chain of inequalities
\begin{center}
$q(G)\leq q(G_1)\leq q(G_3)<q(G_4)\leq q(G_*)=q(K_\delta \vee (K_{n-2\delta +1} \cup (\delta -1)K_1)).$
\end{center}
which contradicts the initial assumption that $q(G)\geq q(K_\delta \vee (K_{n-2\delta +1} \cup (\delta -1)K_1))$.

Since the assumption that $G$ contains no even factor leads to a contradiction in all considered cases, we conclude that $G$ must indeed contain an even factor.
\end{proof}


\vspace{5mm}
\noindent
{\bf Declaration of competing interest}
\vspace{3mm}

The authors declare that they have no known competing financial interests or personal relationships that could have appeared to influence the work reported in this paper.

\vspace{5mm}
\noindent
{\bf Data availability}
\vspace{3mm}

No data was used for the research described in the article.

\vspace{5mm}
\noindent
{\bf Acknowledgements}
\vspace{3mm}

The authors thank genuinely the anonymous reviewers for their helpful comments on this paper.



\end{document}